\documentstyle{amsppt}
\hyphenation{con-sidered}
\nologo
\magnification=\magstep1
	\hcorrection{0.25in}
\leftheadtext{Kenneth A. Ribet}
\rightheadtext{A variant of Fermat's equation}
\newcount\sectCount
\edef\nonouterhead{\noexpand\head}
\def\section#1\endhead{\advance\sectCount by 1%
\nonouterhead\the\sectCount. #1\endhead\noindent}

\newcount\refCount
\def\newref#1 {\advance\refCount by 1
\expandafter\edef\csname#1\endcsname{\the\refCount}}

\edef\nonouterproclaim{\noexpand\proclaim}
\newcount\theoremCount
\def\theorem{\advance\theoremCount
  by 1\nonouterproclaim{Theorem \the\theoremCount}}
\newcount\propCount
\def\proposition{\advance\propCount
  by 1\nonouterproclaim{Proposition \the\propCount}}
\newcount\conjCount
\def\conjecture{\advance\conjCount by 1\definition{Conjecture \the\conjCount}}

\newref ANTWERPIV
\newref DARMONONE
\newref DARMONTWO
\newref DARMONONSERRE  \let\DARMONSERRE\DARMONONSERRE
\newref DARMONGRANVILLE
\newref DENES
\newref DIAMOND
\newref DK
\newref DICKSON
\newref DICKSONTEXT
\newref FREY
\newref IR 
\newref KAMIENNY
\newref EIS
\newref PRIMEDEGREE  
\newref NEWDIRECTIONS
\newref MOMOSE
\newref FERMAT 
\newref DUKE
\newref TW
\newref WASSERMAN
\newref WILES


\let\a\alpha
\def\Q{{\bold Q}}
\def\Z{{\bold Z}}
\def\Qbar{\overline{\Q}}
\def\F{{\bold F}}
\def\Fp{\F_p}
\def\Fell{\F_\ell}
\def\xo(#1){X_0(#1)}
\def\jo(#1){J_0(#1)}
\def\go(#1){\Gamma_0(#1)}
\def\hlinefill{\leaders\hrule height 3pt depth -2.5pt\hfill}
\def\emrule{\thinspace\hbox to 0.75em{\hlinefill}\thinspace}
\def\GalQ{\Gal(\Qbar/\Q)}

\def\GL#1{{\bold{GL}}(2,#1)}

\def\Xsplit(p){X_{\roman{split}}(p)}

\def\startproof{\demo{Proof}}

\def\endproof{{\unskip\nobreak\hfil\penalty50\hskip1em
  \hbox{}\nobreak\hfil$\blacksquare$%
  \parfillskip=0pt \finalhyphendemerits=0
  \par}\enddemo}

\def\Weil{Shimura-Taniyama}

\def\leqno(#1){\eqno{\hbox to .01pt{}\rlap{\rm \hskip -\displaywidth(#1)}}}

\def\newop#1
{\expandafter\def\csname #1\endcsname{\mathop{\roman{#1}}\nolimits}}
\newop ord
\newop Gal 
\newop rad
\newop Aut

\topmatter

\title On the equation $a^p + 2^\a b^p + c^p =0$\endtitle
\author Kenneth A. Ribet\endauthor
\affil University of California, Berkeley\endaffil
\address{Mathematics Department, Evans Hall, University of California,
  Berkeley, CA 94720-3840 USA}\endaddress
\email ribet\@math.berkeley.edu \endemail

\abstract\nofrills
We discuss the equation $a^p + 2^\a b^p + c^p =0$ in which $a$,
$b$, and $c$ are non-zero relatively prime
integers, $p$ is an odd prime number, and $\a$ is a positive
integer.  The technique used to prove Fermat's Last Theorem
shows
that the equation has no solutions
with $\a>1$ or $b$ even.
When $\a=1$ and $b$
is odd, there are the two trivial solutions $(\pm 1, \mp 1, \pm 1)$.
In 1952, D\'enes conjectured that these are the only ones.
Using methods of Darmon, we
prove this conjecture for $p\equiv1$ mod~4.
We link the case
$p\equiv3$ mod~4 to conjectures of Frey and Darmon
about elliptic curves over~$\Q$ with isomorphic
mod~$p$ Galois representations.
\endabstract
\thanks
This article was
prepared while the author was a
research professor at the MSRI, where
research is
supported in part by NSF grant DMS-9022140.
This work was
further
supported by the investigator's
NSF Grant DMS 93-06898.
It is a pleasure to thank 
R.~Borcherds,
H.~Darmon
S.~Kamienny,
I.~Kaplansky,
B.~Mazur and
R.~Tijdeman for helpful feedback and information.
\endthanks
\endtopmatter

\document
\section Introduction\endhead
Let $p\ge5$ be a prime number.
One knows
that Fermat's equation $a^p + b^p + c^p = 0$
has no non-zero integral solutions.
Indeed, suppose that $a^p + b^p + c^p = 0$,
where $a$, $b$ and~$c$ are non-zero.
Following
G.~Frey, one
considers the elliptic curve $E$ with equation
$y^2=x(x-a^p)(x+b^p)$.
The curve $E$ is simultaneously modular~\cite{\WILES, \TW}
and non-modular~\cite{\FERMAT}.
Therefore
no triple $(a,b,c)$ with the hypothesized properties
could have existed.

Ever since A.~Wiles's 1993 announcement that 
Fermat's Last Theorem can be proved along these
lines,
it has been
clear that the proof sketched above can be adapted
to other Diophantine equations having
the skeletal form $A+B=C$.
In particular, 
suppose that $L$ is a prime number taken from the set
$$\Sigma= \{ \, 3, 5, 7, 11, 13, 17, 19, 23, 29, 53, 59 \, \}.$$
The analysis of
J.-P.~Serre~\cite{\DUKE, \S4.3}, combined with
the author's theorem~\cite\FERMAT\
and the
recent work of
Wiles~\cite\WILES\ and Taylor-Wiles~\cite\TW,
provides information
about the family of equations
$ a^p + L^\a b^p + c^p = 0 $.
\theorem
Suppose that $p$ and~$L$
are distinct
prime numbers, with 
$p\ge11$
and
$L\in\Sigma$.
If $\a\ge0$, then there are no triples of non-zero integers
$(a,b,c)$
which satisfy $ a^p + L^\a b^p + c^p = 0 $.
\endproclaim
\noindent
The proof of this theorem can again be summarized succinctly.
A non-zero solution to $ a^p + L^\a b^p + c^p = 0 $
would define a semistable elliptic curve $E$; this
curve would be modular
by~\cite{\WILES, \TW}.
The group of $p$-division points on~$E$ would define
an irreducible two-dimensional representation $\rho$ of~$\GalQ$
over~$\Fp$ with very limited ramification.
Moreover, $\rho$ would be modular because $E$ is modular.
An application of the main theorem
of~\cite\FERMAT\
would lead to the statement that $\rho$ arises from the space
of weight-two cusp forms on~$\go(2L)$.  
As Serre explains in~\cite\DUKE, one may deduce a
contradiction from this statement for
$L$ in~$\Sigma$.

This article concerns the case $L=2$, i.e., the
equation
$$  a^p + 2^\a b^p + c^p =0\leqno($\star$)$$
when $p$ is an odd prime.  
\footnote{Our study of this equation was suggested
by a letter from J.~W.~Weidenman concerning
Diophantine equations which are special cases of ($\star$).
The author wishes to thank him for this inquiry.}
This equation is qualitatively different from those
considered by~Serre,
since ($\star$) has the non-zero solutions
$a=c=-b$ with~$\a=1$.
Their presence is connected up with the fact that the
elliptic curves $E$ defined by solutions to~($\star$)
are not necessarily semistable.
In order to proceed with our analysis,
we must exploit
the fact that
the \Weil\
conjecture holds for all elliptic curves
over~$\Q$ defined by equations of the form $y^2=x(x-A)(x+B)$,
even those curves which are not semistable.
As K.~Rubin and A.~Silverberg have observed,
this extension of Wiles's theorem
follows easily from F.~Diamond's
refinement~\cite\DIAMOND\ of the work of Wiles and Taylor-Wiles.
Alternatively, a somewhat simplified proof 
of the extended theorem
has been given by Diamond and~Kramer~\cite\DK;
these authors appeal directly to~\cite{\WILES, \TW}, rather than
to~\cite\DIAMOND.

In our analysis,
we take
$\a$
to be an integer between 1 and $p{-}1$
without loss of generality.
Also, for technical reasons
we exclude the case $p=3$.
The reader interested in this
omitted
case may consult
Vol.~II, pp.~572--573 of
Dickson's {\it History of the Theory of Numbers}~\cite{\DICKSON}.
According to this {\it History}, Euler showed that $a^3+4b^3+c^3=0$
has no solutions in non-zero integers, while
Legendre established that
$a^3+2b^3+c^3=0$
has only the trivial solutions with $a=c=-b$.
Recently, H.~Wasserman \cite\WASSERMAN\
has communicated a proof of Euler's
result which is inspired by the treatment of
$a^3+b^3+c^3=0$ given by Ireland and Rosen
in~\cite{\IR, Ch.~17, \S8}.

When $\a=1$,
the equation $a^p + 2b^p + c^p=0$
states that the three perfect
$p$th powers
$a^p$, $(-b)^p$ and $c^p$ form an
arithmetic progression.  As the author
learned from
R.~Tijdeman,
there has been considerable interest in arithmetic progressions
consisting of perfect $n$th powers.
While it is easy to exhibit three perfect squares
which form an arithmetic progression (e.g., $7^2$, $13^2$ and~$17^2$),
Fermat stated and
Euler (among others) proved that four distinct squares cannot
form an arithmetic progression. (For a discussion, see
\cite{\DICKSON, Vol.~II, Ch.~XIV}.)
Furthermore, as
Dickson reports in~\cite{\DICKSON, Vol.~II, Ch.~XXII},
Euler proved that $2a^4\pm 2b^4$ is a perfect square
only when $a=b$; in particular,
three distinct fourth powers cannot form
an arithmetic progression.
(For a proof of this latter fact, cf.~\cite{\DICKSONTEXT, Ex.~4, p.~43}.)

For $p$th powers (where $p$ is an odd prime),
D\'enes~\cite\DENES\ made the following conjecture
in~1952:
\conjecture
\edef\conjectureone{\the\conjCount}
Let $p$ be an odd prime.  If $x$, $y$ and~$z$ are
non-zero integers such that $x^p$, $y^p$ and~$z^p$
form an arithmetic progression, then $x$, $y$ and~$z$
are all equal.
\enddefinition
\noindent
Conjecture~\conjectureone\
amounts to the statement that the
only solutions to $a^p+2b^p+c^p=0$ in non-zero integers
are those for which $a=-b=c$.
In support of the conjecture, D\'enes proved the following
theorem~\cite{\DENES, Satz~9}, which implies the conjecture
for all odd primes $p<31$.

\theorem
Suppose that $p$ is a regular odd prime for which the order
of~$2$ in $(\Z/p\Z)^\ast$ is either an even number or else
equal to $p-1\over2$.  Suppose further that $2^{p-1}\not\equiv1$
modulo~$p^2$.  Then the conjecture is true for~$p$.
\endproclaim
\noindent
We prove two theorems about the family ($\star$):

\theorem
\edef\labelA{\the\theoremCount}
The equation $a^p+2^\a b^p+c^p=0$ has no solution in non-zero
integers $a$, $b$, $c$ if $\a$ satisfies $2\le\a<p$.
Furthermore, there are no solutions to
$a^p+2b^p+c^p=0$
in relatively prime
non-zero integers for which $2$ divides $abc$.\endproclaim
\noindent
Given that all elliptic curves $y^2=x(x-A)(x+B)$
are modular, we obtain
Theorem~\labelA\
by mimicing the proof of Fermat's Last
Theorem which we sketched above.
\theorem
\edef\labelB{\the\theoremCount}
If $p\equiv1$ \rom{mod}~$4$, then Conjecture~\conjectureone\
is true for~$p$.
\endproclaim
\noindent
Theorem~\labelB\ is proved by techniques introduced
by Darmon
in~\cite{\DARMONONE, \DARMONTWO}. 
(See also the discussions in~\cite{\DARMONONSERRE, \S4}
and~\cite{\DARMONGRANVILLE, \S4.3}.)
The condition
$p\equiv1$ mod~$4$ in Theorem~\labelB\ is needed so that we
can apply the work of B.~Mazur~\cite\PRIMEDEGREE,
F.~Momose~\cite\MOMOSE, and 
S.~Kamienny~\cite\KAMIENNY\ on the rational points
of modular curves associated with split 
Cartan subgroups of $\GL{\Fp}$.
In fact, we require 
as well
the secondary hypothesis $p\ge17$ to apply
this work, so we do not prove Theorem~\labelB\
for $p=5$ or~$p=13$.  Fortunately, these two primes
are covered by
D\'enes's work.

It is perhaps worth stressing that the hypothesis $p\equiv1$
mod~4 will disappear as soon as theorems for
non-split Cartan subgroups become available.
\section Frey curves\endhead
Let $p$ be an odd prime number.
We view ($\star$) as an equation in the three variables $a$,
$b$ and~$c$ with an auxiliary parameter, $\a$.
We can and do assume that we have
$0<\a<p$.
Suppose that $(a,b,c)$ is a solution to~($\star$) in non-zero
relatively prime integers. It is immediate then that
$a$ and~$c$ are odd; i.e., the three monomials
$A=a^p$, $B=2^\a b^p$ and~$C=c^p$ are relatively prime.
Thus, the
congruence $a\equiv -1$ mod~4 will be satisfied after
possibly multiplying $(a,b,c)$ by~$-1$.
We shall
normalize our solutions by
imposing this congruence.
With this normalization in place, the trivial solutions
$a=c=-b$ with $\a=1$
are reduced to the single triple
$(a,b,c)=(-1,1,-1)$.

Given a normalized solution of~($\star$), one forms
the Frey elliptic curve $E$ with equation
$$ y^2 = x(x-A)(x+B). $$
\theorem
The elliptic curve $E$ is modular.
\endproclaim\noindent
As indicated above,
this theorem was pointed out by
Rubin and Silverberg, who deduced it as
a consequence of the results of~\cite\DIAMOND.
After learning of the Rubin-Silverberg observation,
F.~Diamond and K.~Kramer gave a more ``elementary'' proof
of the theorem in~\cite\DK.  This latter article applies the
work of Wiles and Taylor-Wiles, but does not rely on the
refinements of~\cite\DIAMOND.
It contains a great deal of
information about the arithmetic of Frey curves,
some of which we shall
recall below.

Because of our normalization, the integer $A$ satisfies
$A\equiv -1$ mod~4; furthermore,
$B$ is even.
These are the conventions that were employed
in~\cite\DUKE\ and~\cite\DK.
The calculations of~\cite{\DUKE, \S4.1}
show
that the conductor $N_E$
of~$E$ has the form $2^t\rad'(ABC)$ where $t$ is a non-negative
integer.  Here, 
we have written
$\rad'(ABC)$ for
the product of the odd prime
divisors of~$ABC$.  In particular,
the curve $E$ is semistable at all primes $p\ne2$.
The precise value of~$t$ is computed by Diamond and
Kramer~\cite\DK, who find
that $t$ is
5, 3, 3, 0 or~1 according as $\ord_2(B)$ is 1, 2, 3, 4,
or an integer greater than~4.  Thus $E$ is
semistable at~2 if and only if $B$ is divisible by~16.
Since $E$ is in any case semistable away from~2,
$E$ is a semistable elliptic curve precisely
when 16 divides~$B$.
The
minimal discriminant $\Delta_E$ of~$E$ may be written
$2^u(ABC)^2$ where $u$ is an integer which is calculated
in~\cite\DK.
For instance,
$u=-8$ when $t=1$.
Therefore
$$ \ord_\ell(\Delta_E)\equiv0 \hbox{ mod }p$$
for all primes $\ell\ne2$.

\proclaim{Lemma}
\edef\labelS{\the\sectCount}
The conductor $N_E$ is a power of~2 if and only if $(a,b,c)$
is the trivial solution $(-1,1,-1)$.
\endproclaim
\startproof
The solution $(-1,1,-1)$ to~($\star$) for the value $\a=1$
leads to the elliptic curve $E=E_0$ with equation
$y^2 = x(x+1)(x+2)$. 
A translation in~$x$
transforms $E_0$
into
the familiar complex multiplication elliptic curve
$y^2 = x^3-x$ of conductor~$32$.
Conversely, suppose that $N_E$ is a power of~2.  Then
$\rad'(ABC)=1$, so that $ABC$ is a power of~2.  Since $a\equiv-1$
mod~4, we have $a=-1$.  Similarly, the odd number
$c$ can only be $\pm1$ and
$b$ must be a power of~2.  The equation
$-1 + 2^\a b^p + (\pm1) =0$ forces $b=1$, $\a=1$ and $\pm1 = -1$.
\endproof
\proclaim{Corollary}
If $(a,b,c)$ is not the trivial solution, then $E$ has
multiplicative reduction at some prime $q\ne2$.
\endproclaim
\startproof
This is clear since $N_E$ is a power of~2 times a square-free
odd number.
\endproof

For each prime number $\ell$, let $E[\ell]$ be the group of
$\ell$-division points on~$E$, regarded as a two-dimensional
representation of~$\GalQ$ over the field $\Fell$.  We recall
the following fact.

\proposition
The representation $E[\ell]$ is irreducible for all primes $\ell\ge5$.
Moreover, if $E$ is not semistable, then $E[3]$ is irreducible.
\edef\labelC{\the\propCount}
\endproclaim

\startproof
First suppose that $E$ is semistable over~$\Q$.  Then,
as was noted in~\cite\DUKE, 
the
result to be proved follows easily from a theorem
of Mazur~\cite{\EIS, \PRIMEDEGREE}.
More precisely,
suppose that $\ell\ge5$ and that $E[\ell]$ is
reducible.  Then $E$ has a rational subgroup $C$ of order~$\ell$.
The semistability hypothesis implies that the action of~$\GalQ$
on~$C$ is ramified only at~$\ell$, and a local study at~$\ell$
then shows that $\GalQ$ must act on~$C$ either trivially or
via the mod~$\ell$ cyclotomic character.  This implies that
some elliptic curve over~$\Q$ which is isogenous to~$E$ contains a
group of rational points which is isomorphic to $\Z/2\Z\oplus
\Z/2\ell\Z$.  The existence of such a curve is incompatible with
Ogg's Conjecture, which was proved by Mazur in~\cite\EIS.

Now suppose that $E$ is not semistable; this means that $E$
has additive reduction at~2.  Then the indicated irreducibility
follows from a stronger statement which is proved
by Diamond and Kramer in~\cite\DK:
Let $I$ be an inertia subgroup of~$\GalQ$ for the prime~2;
then the action of~$I$ on~$E[\ell]$ is irreducible
if $\ell\ge3$.  Since the proof
of this statement
is quite elementary, we shall
recall it now for the convenience of the reader.  

Since $E$ has additive reduction at~2, the 2-part of $N_E$ may
be written $2^{2+\delta}$, where $\delta$ is the 
exponent of the Swan conductor
of the representation given by the action of~$I$ on~$E[\ell]$.
As we noted above,
$2+\delta$ is equal to either 5 or~3; thus $\delta$ is an odd number.
Assume now that $E[\ell]$ is reducible as an $I$-module.
Then $E[\ell]$ is an extension of one 1-dimensional representation
by another, and $\delta$ is the sum of the conductors of the
two characters associated with the 1-dimensional
representations.  These characters are in fact inverses of each
other, since $I$ acts trivially on the determinant of~$E[\ell]$.
(This determinant corresponds to the mod~$\ell$ cyclotomic
character, which is unramified at~2.)  Hence the conductors of
the two characters are equal, giving that $\delta$ is even.
\endproof
\proclaim{Corollary}
Suppose that $p\ge5$ or that $p=3$ and $b$ is odd.
Then $E[p]$ is irreducible.
\endproclaim
\startproof
The only point to be checked is that $E$ is non-semistable
if $p=3$ and $b$ is odd.  In fact, suppose that $p=3$.
Then $E$ is semistable if and only if $b$ is even.  Indeed,
if $b$ is even, then 8 divides $b^p$, so that 16 divides $B$.
Conversely, suppose that 16 divides $B=2^\a b^p$.  Since $1\le\a\le2$,
it is clear that $b$ is even.
\endproof

\section Proofs of Theorem \labelA\ and \labelB\endhead
Suppose that
$a^p+2^\a b^p+c^p=0$,
where the integers $a$, $b$ and $c$ are non-zero and relatively
prime, and where $\a$ satisfies $1\le\a<p$.
It is evident then that $a$ and $c$ are odd.
As above, we multiply $(a,b,c)$ by~$-1$ if necessary
in order to ensure that $a$ is congruent to~$3$ mod~$4$.
We again form the Frey curve
$E: y^2=x(x-A)(x+B)$. In the notation introduced above
the conductor $N_E$ of~$E$ is the product $2^t\rad'(ABC)$,
for some integer $t$ 
in the set $\{\, 0 , 1, 3, 5\,\}$.
We have $t\le3$ if and only if the even number
$B$ is divisible by~$4$; we have $t=5$ in the contrary case.
We will suppose from now on that $p\ge5$.

We first prove Theorem~\labelA, i.e., that $t=5$.
Because $p\ge5$,
we may
deduce from the Corollary
to Proposition~\labelC\
that the representation
$$ \rho\:\GalQ\to\GL{\Fp}$$
defined by~$E[p]$ is irreducible.
It is modular of level~$N_E$ (i.e., it arises from the
space of weight-two cusp forms on~$\go(N_E)$) because
$E$ is a modular elliptic curve of conductor~$N_E$.
Since $\Delta_E$ is a perfect $p$th power times a power of~2, the
representation $\rho$ is finite at each prime $\ell\ne2$.  
The main theorem of~\cite\FERMAT\ thus implies that $\rho$ is
modular of level~$2^t$.  (Each odd prime $\ell$ dividing~$N_E$
can be jettisoned from the level of~$\rho$.)
We  conclude that $t=5$, since
there are no non-zero cusp forms of weight~two
on~$\go(8)$.  Equivalently, $\ord_2(B)=2$, 
as asserted by Theorem~A.

\remark{Remark}
In the omitted case $p=3$, suppose that $b$ is
odd.  Then $E[3]$ is again an irreducible representation,
and the argument we have given may be used to deduce
that 
$\a=1$.
\endremark

Continuing the argument, we now prove Theorem~\labelB.
Let
$E_0$ again be the elliptic
curve which is associated with the trivial solution $(-1,1,-1)$,
i.e., the elliptic curve over~$\Q$ with equation $y^2=x^3-x$.
\proposition
The two-dimensional mod~$p$ representations of~$\GalQ$
which are defined by $E$ and~$E_0$
are isomorphic.
\edef\labelD{\the\propCount}
\endproclaim
\startproof
Let $\rho$ be the mod~$p$ representation which is defined
by~$E$, i.e., by the space $E[p]$ of $p$-division points of~$E$;
let $\rho_0$ be the analogue of~$\rho$ for~$E_0$.
We have seen that the irreducible representation
$\rho$ is associated with an eigenform in the space
of weight-two cusp forms on~$\go(2^t)=\go(32)$.
It is a known fact that this space
is one-dimensional; equivalently, $\jo(32)$
is an elliptic curve.
(See, e.g., \cite{\ANTWERPIV, p.~136}.)
It follows that $\rho$ is the mod~$p$
representation 
$\jo(32)[p]$.
In particular, 
the isomorphism class of~%
$\rho$ is independent of the solution
$(a,b,c)$ giving rise to~$E$.
Therefore, $\rho$ and~$\rho_0$ are isomorphic, as stated.
\endproof

We next recall
the well known fact that
the image of~$\rho_0$ is contained in the normalizer of
a Cartan subgroup of~$\GL{\Fp}$.
Indeed, let $R=\Z[\mu_4]$ be the full ring of endomorphisms
of~$E_0$.  Then $E_0[p]$ is easily seen to be a free rank-1
module over $R/pR$.  Let $C$ be the image of $(R/pR)^*$
in the group of automorphisms of~$E_0[p]$, so that $C$ is
either $\F_p^*\times\F_p^*$ or~$\F_{p^2}^*$, according
as $p$ is congruent to~1 or to~$-1$ mod~4.  Then $C$ is
a Cartan subgroup of~$\Aut E_0[p]\approx \GL{\Fp}$.
One says that $C$ is split or
non-split according as $p$ is 1 or$-1$ mod~4.
The restriction of~$\rho_0$
to~$\Gal\bigl(\Qbar/\Q(\sqrt{-1})\bigr)$ takes values in~$C$,
and the full image of~$\rho$ takes values in the normalizer
of~$C$ in~$\GL{\Fp}$.  The index of $C$ in its normalizer is~2.

Suppose now that $p\equiv1$ mod~4.
Then $E$ defines a point on the modular curve denoted $\Xsplit(p)$,
cf.~\cite{\EIS, Ch.~III, \S6} and the discussions in \S3 and~\S4e
of~\cite\PRIMEDEGREE.  This circumstance puts strong constraints
on the set of prime numbers dividing the 
denominator of the
$j$-invariant of~$E$,
i.e., the set of primes at which $E$ does not have potential
good reduction.  Specifically,
if $p\ge17$,
a result of Mazur \cite{\PRIMEDEGREE, Cor.~4.8}
proves that $E$ has potential good reduction at all primes $\ell\not=2,p$
satisfying $\ell\not\equiv\pm1$ mod~$p$.
(This result holds also for $p=11$, but this is irrelevant to
our application, which requires $p\equiv1$ mod~4.)
Mazur's theorem has been strengthened by subsequent work.
In particular, F.~Momose \cite{\MOMOSE, Prop.~3.1}
proves that $E$ has potential good reduction at all primes $\ell\ne2$,
as long as the prime $p$ satisfies $p\ge17$.

Alternatively, 
under the same hypothesis on~$p$,
Darmon notes in~\cite{\DARMONONE, Cor.~1.7}
that $E$ has potential good reduction at all primes $\ell\ne2,3$;
this observation is obtained
by combining
a theorem
of Kamienny~\cite\KAMIENNY\ with~%
\cite{\PRIMEDEGREE, Cor.~4.3}.
Darmon's result 
concerns elliptic curves over~$\Q(\sqrt{-1})$ and
requires only that $E$ possess a rational subgroup
of order~$2p$ over this field.

Suppose now that we have $p\ge17$ and $p\equiv1$ mod~4.
Then if $(a,b,c)$ is a normalized solution to $a^p+2b^p+c^p=0$,
the corresponding curve $E$ has multiplicative reduction
at all odd primes $\ell$ dividing $abc$.  
Momose's result implies that there is no such prime; Darmon's
implies that 3 is the only possible such prime.
On either count,
we find that 
two of
$a$, $b$ and~$c$ 
are $\pm1$ while the third is $\pm3^n$ for some $n\ge0$.
Indeed, $a$, $b$ and~$c$ are relatively prime and all
of them are odd in view of Theorem~\labelA.
Elementary reasoning allows us to reach a contradiction.

\section A conjecture of Frey\endhead
\conjecture\edef\junk{\the\conjCount}
Let $A$ be an elliptic curve over~$\Q$.
Then all sufficiently large prime numbers $p$
have the following property:
if $B$ is an elliptic curve over~$\Q$
for which $A[p]$ and $B[p]$ are
isomorphic representations
of~$\GalQ$, then $A$ and~$B$ are isogenous over~$\Q$.
\enddefinition
\noindent
Conjecture~\junk\ 
appears as Conjecture~4.3 in~\cite\DARMONSERRE, where
it is attributed to G.~Frey.
It is similar in flavor to the conjectural statements
in Frey's article~\cite\FREY.
The reader is invited to consult Mazur's article~\cite\NEWDIRECTIONS\
as well as \cite\DARMONONSERRE\ and~\cite\FREY\ for variants
and generalizations.  Here is one such generalization~\cite{\DARMONSERRE,
Conj.\ 4.4 and Conj.~4.5}:

\conjecture
There is an integer $t>0$
with the following property.
Suppose that $A$
and $B$ are elliptic curves over~$\Q$
and that the
Galois representations $A[p]$ and~$B[p]$
have isomorphic semisimplifications.
If $p>t$, then
$A$ and $B$ are isogenous.
\enddefinition

We record the following simple observation:
\proposition
Suppose that Conjecture~\junk\
is true.
Then Conjecture~\conjectureone\
holds for all
sufficiently large prime numbers~$p$.
\endproclaim

\startproof
Suppose that $(a,b,c)$ is a normalized solution to
$a^p+2b^p+c^p=0$.  If $E$ is the associated Frey curve,
then we have seen that $E[p]$ and~$E_0[p]$ are isomorphic.
Applying Conjecture~\junk\ with $A=E_0$, we find that $E$
and~$E_0$ are isogenous for $p$ sufficiently large.
The isogeny relation between $E$ and~$E_0$ implies that
these two elliptic curves have the same primes of bad reduction,
so that $E$ has good reduction outside~2.  By the lemma of~\S\labelS,
this implies that $(a,b,c)=(-1,1,-1)$.
Hence $a^p+2b^p+c^p=0$ has only the trivial normalized solution
for sufficiently large~$p$.
\endproof

\enddocument